\documentclass[10pt,a4paper]{article}
\usepackage[ascii]{inputenc}
\usepackage[T1]{fontenc}
\usepackage{amsmath,amssymb,amsfonts,amsthm}
\usepackage[a4paper,bottom=3.4cm]{geometry}
\usepackage{graphicx}
\usepackage{authblk}
\usepackage[caption=false]{subfig}
\usepackage[pdftex,bookmarks=false,colorlinks=true,linkcolor=blue,
citecolor=blue,filecolor=black,urlcolor=blue]{hyperref}

\DeclareMathOperator{\e}{e}
\DeclareMathOperator{\Tr}{Tr}
\DeclareMathOperator{\erf}{erf}

\newtheorem{theorem}{Theorem}
\newtheorem*{remark}{Remark}

\newcommand{\ud}{\mathrm{d}}
\newcommand{\abs}[1]{\lvert{#1}\rvert}
\newcommand{\norm}[1]{\lVert{#1}\rVert}
\newcommand{\vect}[1]{\mathbf{#1}}

\newcommand{\N}{\mathbb{N}}
\newcommand{\Z}{\mathbb{Z}}
\newcommand{\R}{\mathbb{R}}
\newcommand{\C}{\mathbb{C}}

\begin{document}

\title{Efficient numerical evaluation of thermodynamic quantities on infinite (semi-)classical chains}

\author[1]{Christian B.~Mendl\footnote{\href{mailto:christian.mendl@tum.de}{christian.mendl@tum.de}}}
\author[2]{Folkmar Bornemann\footnote{\href{mailto:bornemann@tum.de}{bornemann@tum.de}}}
\affil[1]{\normalsize Technical University of Munich, Department of Informatics and Institute for Advanced Study, Boltzmannstr.~3, 85748 Garching, Germany}
\affil[2]{\normalsize Technical University of Munich, Department of Mathematics, Boltzmannstr.~3, 85748 Garching, Germany}

\date{\normalsize June 24, 2020}

\maketitle

\begin{abstract}
This work presents an efficient numerical method to evaluate the free energy density and associated thermodynamic quantities of (quasi) one-dimensional classical systems, by combining the transfer operator approach with a numerical discretization of integral kernels using quadrature rules. For analytic kernels, the technique exhibits exponential convergence in the number of quadrature points. As demonstration, we apply the method to a classical particle chain, to the semiclassical nonlinear Schr\"odinger equation and to a classical system on a cylindrical lattice.
\end{abstract}

\section{Introduction}

The partition function $Z$ is a cornerstone of statistical mechanics \cite{LandauLifshitz1980, Huang2010}, in particular since thermodynamics quantities can be obtained as derivatives of $Z$. Nevertheless, computing $Z$ for larger systems is a challenging task if a closed-form solution cannot be found. In this work we are concerned with an efficient numerical method for evaluating $Z$ and the associated free energy density for (semi-)classical chains, in the thermodynamic limit of infinite system size. In more detail, we consider a Hamiltonian $H$ on a (quasi) one-dimensional lattice of size $L$, containing solely nearest neighbor terms:
\begin{equation}
\label{eq:formal_hamiltonian}
H(z^1, \dots, z^L) = \sum_{\ell=1}^L h(z^\ell, z^{\ell+1}),
\end{equation}
with $z^\ell \in \Omega \subset \R^n$ for all $\ell$ and some fixed $n \in \N$. We assume periodic boundary conditions and thus set $z^{L+1} = z^1$. A canonical example is a classical particle chain, with $z^\ell = (p_\ell, q_\ell) \in \R^{2d}$ the momentum and position of the $\ell$-th particle ($d \in \{1, 2, 3\}$), and
\begin{equation}
\label{eq:Hpc}
H_{\text{pc}}(p_1, \dots, p_L, q_1, \dots, q_L) = \sum_{\ell=1}^L \left(\tfrac{1}{2} p_\ell^2 + V(q_\ell, q_{\ell+1}) \right)
\end{equation}
consisting of kinetic energy terms $\tfrac{1}{2} p_\ell^2$ and the potential energy described by $V$.

The partition function of the physical system is defined as
\begin{equation}
Z_L(\beta) = \int_{\Omega} \cdots \int_{\Omega} \e^{-\beta H(z^1, \dots, z^L)} \ud z^1 \cdots \ud z^L,
\end{equation}
with $\beta = 1/(k_{\text{B}} T) \in \R^{+}$ the ``inverse temperature''. Our goal is to take the thermodynamic limit $L \to \infty$ and compute the canonical free energy density
\begin{equation}
\label{eq:free_energy_def}
F(\beta) = -\beta^{-1} \lim_{L \to \infty} \frac{1}{L} \log Z_L(\beta).
\end{equation}

\section{Method}
\label{sec:method}

To evaluate the partition function numerically, we combine the well-known transfer-matrix method \cite{KramersWannier1941, Mussardo} with a numerical discretization of integral kernels based on numerical quadrature methods. The latter idea traces back to Nystr\"om \cite{Nystrom1930} and was recently applied by the second author to the computation of Fredholm determinants \cite{Bornemann2010}.

\subsection{Assumptions}
\label{subsect:assumpt}

To keep technical preliminaries as simple as possible, we make the following assumptions, which cover the applications in the next section: Let $\Omega$ be a (not necessarily finite) interval of $\R^n$ and $\nu$ be a finite positive measure on $\Omega$ with a strictly positive, continuous density $\omega$ (weight function). Now we rewrite the partition function in the form
\begin{equation}
\label{eq:Z_kbeta_kernel}
Z_L(\beta) = \int_{\Omega} \cdots \int_{\Omega} \prod_{\ell=1}^L k_{\beta}(z^\ell, z^{\ell+1}) \, \ud\nu(z^1) \cdots \ud\nu(z^L)
\end{equation}
with $z^{L+1} = z^1$ as before. Here the integral kernel
\begin{equation}
k_{\beta}: \Omega \times \Omega \to \R, \quad k_{\beta}(z, z') = \frac{\e^{-\beta h(z, z')}}{\sqrt{\omega(z) \omega(z')}}
\end{equation}
is assumed to be continuous and bounded; it is strictly positive by construction. We further assume that $h(z, z') = h(z', z)$, i.e., $k_\beta$ is a symmetric kernel. These assumptions imply, in particular,
\begin{equation}
\int_{\Omega} \int_{\Omega} \lvert k_{\beta}(z, z') \rvert^2 \, \ud\nu(z) \, \ud\nu(z') < \infty,
\end{equation}
so that the symmetric kernel $k_\beta$ induces an integral operator \cite[\S16.1]{Lax}
\begin{equation}
\label{eq:Tbeta_def}
(\mathcal{T}_{\beta} u)(z) = \int_{\Omega} k_{\beta}(z, z') u(z') \, \ud\nu(z'), \quad u \in L^2(\Omega, \nu).
\end{equation}
on the Hilbert space $L^2(\Omega,\nu)$ that is self-adjoint and Hilbert--Schmidt (and hence, compact). Since the product of two Hilbert-Schmidt operators is of trace class and the trace class operators form an ideal within the bounded ones, Eq.~\eqref{eq:Z_kbeta_kernel} for $L \ge 2$ may be represented as\footnote{This holds for $L=1$ only if $\mathcal{T}_{\beta}$ would be trace class itself, cf.\ the discussion in \cite[pp.~878--879]{Bornemann2010}.}
\begin{equation}
\label{eq:Z_trace_Kbeta_op}
Z_L(\beta) = \Tr\!\left(\mathcal{T}_{\beta}^L\right).
\end{equation}

\subsection{Dominant eigenvalue and free energy}

Following the mathematical theory of compact operators on Hilbert spaces \cite[\S21.2]{Lax}, the non-zero elements of the spectrum of $\mathcal{T}_{\beta}$ (which is real since the operator is selfadjoint) are at most countably many eigenvalues of finite multiplicity that accumulate only at $0$ (which belongs to the spectrum).

Generalizing the Perron--Frobenius theory in matrix analysis \cite[\S8.2]{Horn85}, Jentzsch's Theorem \cite[\S8.7, Satz 3]{Fen2} states that the Hilbert-Schmidt operator $\mathcal{T}_\beta$ with \emph{strictly positive} kernel $k_\beta$ has a simple, dominant, strictly positive eigenvalue. That is, all the non-zero eigenvalues can be ordered as (with each eigenvalue listed as often as given by its multiplicity)
\[
\lambda_1(\mathcal{T}_\beta) > \abs{\lambda_2(\mathcal{T}_\beta)} \geq \abs{\lambda_3(\mathcal{T}_\beta)} \geq \cdots > 0.
\]
Hence, evaluating Eq.~\eqref{eq:Z_trace_Kbeta_op} together with Lidskii's theorem \cite[\S30.3]{Lax}, that is,
\[
\Tr\!\left(\mathcal{T}_{\beta}^L\right) = \sum_{j}\lambda_j(\mathcal{T}_{\beta})^L \qquad (L\geq 2),
\]
leads to
\begin{equation}
\lim_{L \to \infty} \frac{1}{L} \log Z_L(\beta) = \lim_{L \to \infty} \frac{1}{L} \log \sum_j \lambda_j(\mathcal{T}_{\beta})^L = \log \lambda_1(\mathcal{T}_{\beta}),\qquad F(\beta) = -\beta^{-1} \log \lambda_1(\mathcal{T}_{\beta}).
\end{equation}
Thus the calculation of the free energy amounts to computing the dominant eigenvalue of $\mathcal{T}_{\beta}$.

\subsection{Nystr\"om method for computing the dominant eigenvalue}
\label{subsect:nystrom}

For simplicity, we restrict the discussion of the numerical method to $n=1$, that is, $\Omega \subset\R$. For the multi-dimensional case $n>1$, see the discussion of the example in Sect.~\ref{subsect:lattice}.

We make use of a numerical quadrature rule of the form
\begin{equation}
\label{eq:quadrature_rule}
\int_\Omega f(z)\,d\nu(z) \approx \sum_{i=1}^m w_i f(z_i)
\end{equation}
with \emph{positive} weights $w_i > 0$ and nodes $z_i$. The approximation is of order $p\geq 1$ if the quadrature rule is exact for polynomials of degree at most $p-1$. It is possible to construct such a unique quadrature rule of maximum order $p=2m$, called Gauss quadrature rule; see \cite{GolubWelsch} for the classical construction of Gaussian weights and nodes from the tridiagonal Jacobi matrix of the orthogonal polynomials associated with the measure $\nu$. Error estimates depend on the regularity of the integrand $f$, e.g., \cite[\S4.8]{Rabi07}: the error for $m\to\infty$ is of the form $O(m^{-k})$ if $f \in C^{k-1,1}$ and of the form $O(e^{-c m})$ with some $c>0$ if $f$ extends analytically to an ellipse or (semi)-strip in the complex plane containing $\Omega$. The former estimate is called \emph{algebraic convergence}, the latter \emph{exponential convergence}.

Inserting the quadrature rule into Eq.~\eqref{eq:Tbeta_def} results in
\begin{equation}
(\mathcal{T}_{\beta} u)(z_i) \approx \sum_{j=1}^m k_{\beta}(z_i, z_j) w_j u(z_j)
\end{equation}
for all $i = 1, \dots, m$. Setting $u_i = \sqrt{w_i} u(z_i)$, we have thus discretized the integral operator by the symmetric $m \times m$ matrix (Nystr\"om method)
\begin{equation}
\label{eq:T_beta_def}
T_{\beta} = \big(k_{\beta}(z_i, z_j)\,\sqrt{w_i\,w_j}\big)_{i,j=1}^m.
\end{equation}
By the Perron--Frobenius theory \cite[Cor.~8.2.6]{Horn85}, this (element-wise) strictly positive matrix has a simple, dominant, strictly positive eigenvalue $\lambda_1(T_\beta)$. The following theorem shows that $\lambda_1(T_{\beta}) \approx \lambda_1(\mathcal{T}_{\beta})$, and thus for the free energy density in Eq.~\eqref{eq:free_energy_def}
\begin{equation}
F(\beta) \approx -\beta^{-1} \log \lambda_1(T_{\beta}),
\end{equation}
with a speed of convergence, algebraic or exponential, depending on the smoothness of $k_\beta$.

\begin{theorem}\label{thm:conv} Under the assumptions in Sect.~\ref{subsect:assumpt} the error $\epsilon_m = \lambda_1(T_{\beta}) - \lambda_1(\mathcal{T}_{\beta})$
of the Nystr\"om method for computing the dominant eigenvalue behaves as follows as $m\to\infty$:
\begin{itemize}
\item if $k_\beta \in C^{k-1,1}$ then there is algebraic convergence $\epsilon_m=O(m^{-k})$;
\item if $k_\beta$ extends analytically to $E\times E$, where $\Omega\subset E\subset\C$ is an ellipse or (semi-)strip, there is exponential convergence  $\epsilon_m=O(e^{-c m})$ with some constant $c>0$.
\end{itemize}
\end{theorem}
The classical, though long and quite involved proof goes by the theory of collectively compact operators and can be found in many books on the numerical treatment of integral equations, e.g., \cite[Thm.~4.8.20]{Hack95}. Note, that though only the case of algebraic convergence is covered in this reference, the proof extends literally to the case of exponential convergence by adjusting the consistence assumptions accordingly.

An alternative, conceptually much simpler novel proof based on the theory of the Fredholm determinant can be found in the Appendix of this paper.

\begin{remark} While the Nystr\"om method for the dominant eigenvalue essentially inherits the convergence properties of an underlying cubature formula also in the multidimensional case $n>1$ (see the proof in the Appendix for guidance), we refrained from formulating a general theorem since the characterization of convergence properties of general cubature formulae is much more involved in the first place. If, however, a Gaussian quadrature is applied coordinate-wise, Theorem~\ref{thm:conv} extends in a straightforward fashion (in the case of exponential convergence $E\subset\C^n$ has then to be chosen as an poly-ellipse or poly-(semi-)strip), cf. the example in Sect.~\ref{subsect:lattice}.
\end{remark}

\section{Applications}

We demonstrate the range of applicability of the method via a diverse selection of model systems.

\subsection{Classical particle chain}
\label{sec:particle_chain}

Consider a classical particle chain governed by a Hamiltonian of the form \eqref{eq:Hpc} above. For evaluating the partition function, the integration over the momentum variables can be performed in closed form, such that
\begin{equation}
Z_L(\beta) = \left(\frac{2\pi}{\beta}\right)^{L/2} \tilde{Z}_L(\beta)
\end{equation}
with
\begin{equation}
\label{eq:particle_chain_tildeZ_def}
\tilde{Z}_L(\beta) = \int_{-\infty}^{\infty} \cdots \int_{-\infty}^{\infty} \prod_{\ell=1}^L \e^{-\beta V(q_\ell, q_{\ell+1})} \ud q_1 \cdots \ud q_L.
\end{equation}
As specific example for the following, we choose
\begin{equation}
\label{eq:particle_chain_V}
V(q, q') = V_{\text{loc}}(q) + \tfrac{1}{2} \gamma (q - q')^2
\end{equation}
with anharmonic on-site potential
\begin{equation}
V_{\text{loc}}(q) = \tfrac{1}{2} \eta q^2 + \tfrac{1}{6} \mu q^3 + \tfrac{1}{24} \lambda q^4
\end{equation}
and coefficients $\lambda, \mu, \eta, \gamma \in \R$, $\eta > 0$, $\abs{\lambda} \ge \abs{\mu}$. The $\sim q^4$ term ensures that $V$ grows asymptotically to infinity as $\abs{q} \to \infty$, and the above $V(q, q')$ is equivalent to its symmetrized version $\frac{1}{2} (V_{\text{loc}}(q) + V_{\text{loc}}(q')) + \frac{1}{2} \gamma (q - q')^2$.

We now assign the term $\frac{1}{2} \eta q^2 $ from the on-site potential to the weight function:
\begin{equation}
\label{eq:particle_chain_weight}
\omega: \R \to \R^+, \quad \omega(q) = \frac{\e^{-\frac{1}{2} \beta \eta q^2}}{\sqrt{2\pi/(\beta \eta)}}.
\end{equation}
This leads to a (rescaled) Gauss-Hermite quadrature rule \cite[\S3.5(v)]{NIST2010},  and we denote the weights by $w_i$ and nodes by $z_i$, $i = 1, \dots, m$, as above.\footnote{A modern, fast numerical method of optimal $O(m)$ complexity for calculating these weights and nodes can be found in \cite{Town16}.} The particular choice of the weight function is somewhat arbitrary --- we could have also included the cubic and quartic terms from the on-site potential into the weight function; at the expense, though, of a less straightforward computation of the weights and nodes of the quadrature rule. The general reasoning is to capture most of the local weight while retaining a well-behaved kernel for the genuine inter-particle potential (see Eq.~\eqref{eq:particle_chain_kernel} below).

Combining the weight function \eqref{eq:particle_chain_weight} with Eq.~\eqref{eq:particle_chain_tildeZ_def} leads to
\begin{equation}
\tilde{Z}_L(\beta) = \left(\frac{2\pi}{\beta \eta}\right)^{L/2} \int_{-\infty}^{\infty} \cdots \int_{-\infty}^{\infty} \prod_{\ell=1}^L k_{\beta}(q_\ell, q_{\ell+1}) \, \ud\nu(q_1) \cdots \ud\nu(q_L)
\end{equation}
with the symmetrized kernel
\begin{equation}
\label{eq:particle_chain_kernel}
k_{\beta}(q, q') = \e^{-\frac{1}{12} \beta \mu (q^3 + {q'}^3) - \frac{1}{48} \beta \lambda (q^4 + {q'}^4) - \frac{1}{2} \beta \gamma (q - q')^2}.
\end{equation}
We then assemble the symmetric matrix in Eq.~\eqref{eq:T_beta_def}. Finally, after taking all prefactors in to account, the numerical approximation of the free energy density reads
\begin{equation}
- \beta F(\beta) = \lim_{L \to \infty} \frac{1}{L} \log Z_L(\beta) \approx \log\frac{2\pi}{\beta} - \frac{1}{2} \log\eta + \log \lambda_1(T_{\beta}).
\end{equation}

The special case $\gamma = 0$ serves as reference, since the partition function factorizes in this case, i.e., $Z_L(\beta)\vert_{\gamma = 0} = (Z_1(\beta)\vert_{\gamma = 0})^L$ with
\begin{equation}
\label{eq:particle_chain_Z1_gamma0}
\log Z_1(\beta)\vert_{\gamma = 0} = \frac{1}{2} \log \frac{2\pi}{\beta} + \log \int_{-\infty}^{\infty} \e^{-\beta V_{\text{loc}}(q)} \ud q.
\end{equation}

\begin{figure}[!tb]
\centering
\subfloat[$\gamma = 0$]{\includegraphics[width=0.48\textwidth]{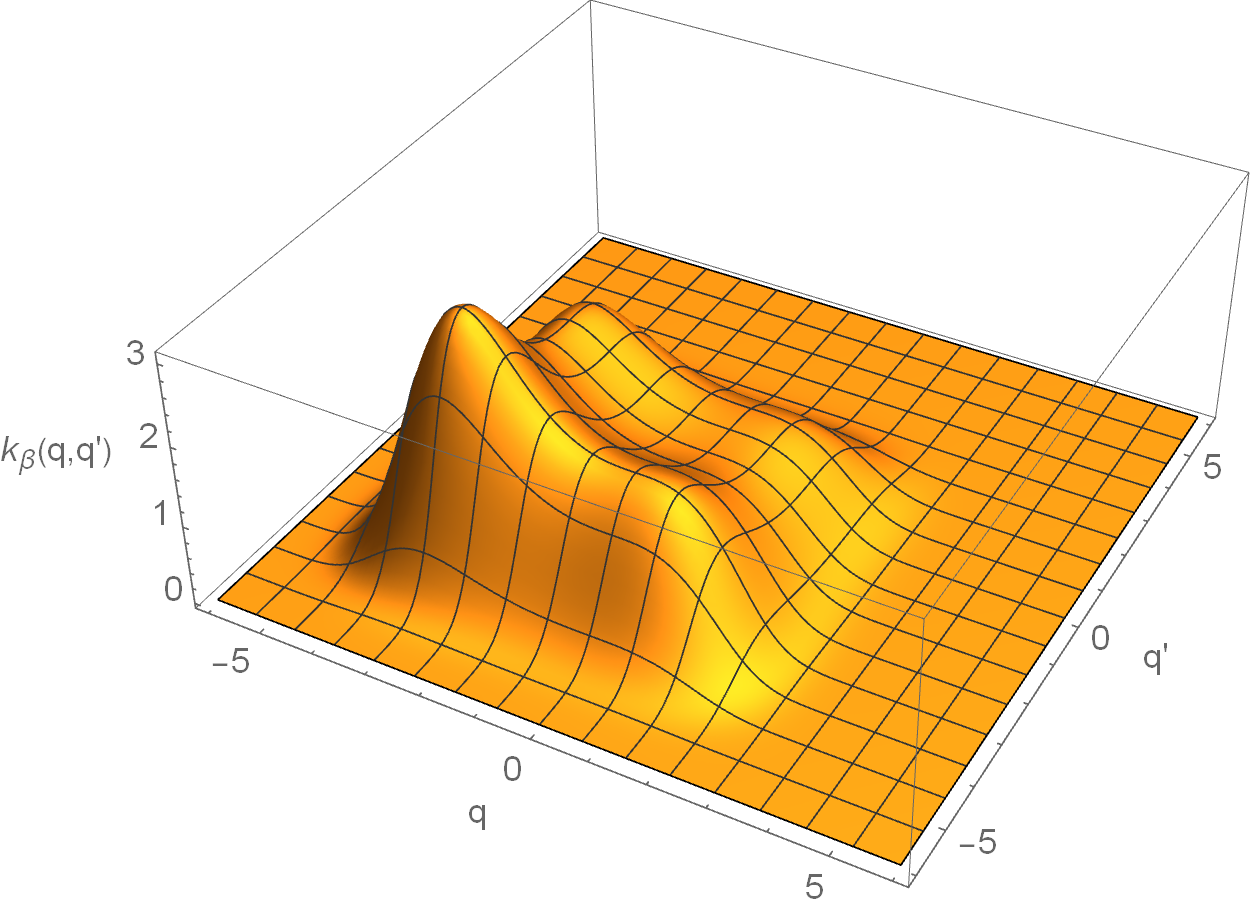}}%
\hspace{0.01\textwidth}%
\subfloat[$\gamma = 1$]{\includegraphics[width=0.48\textwidth]{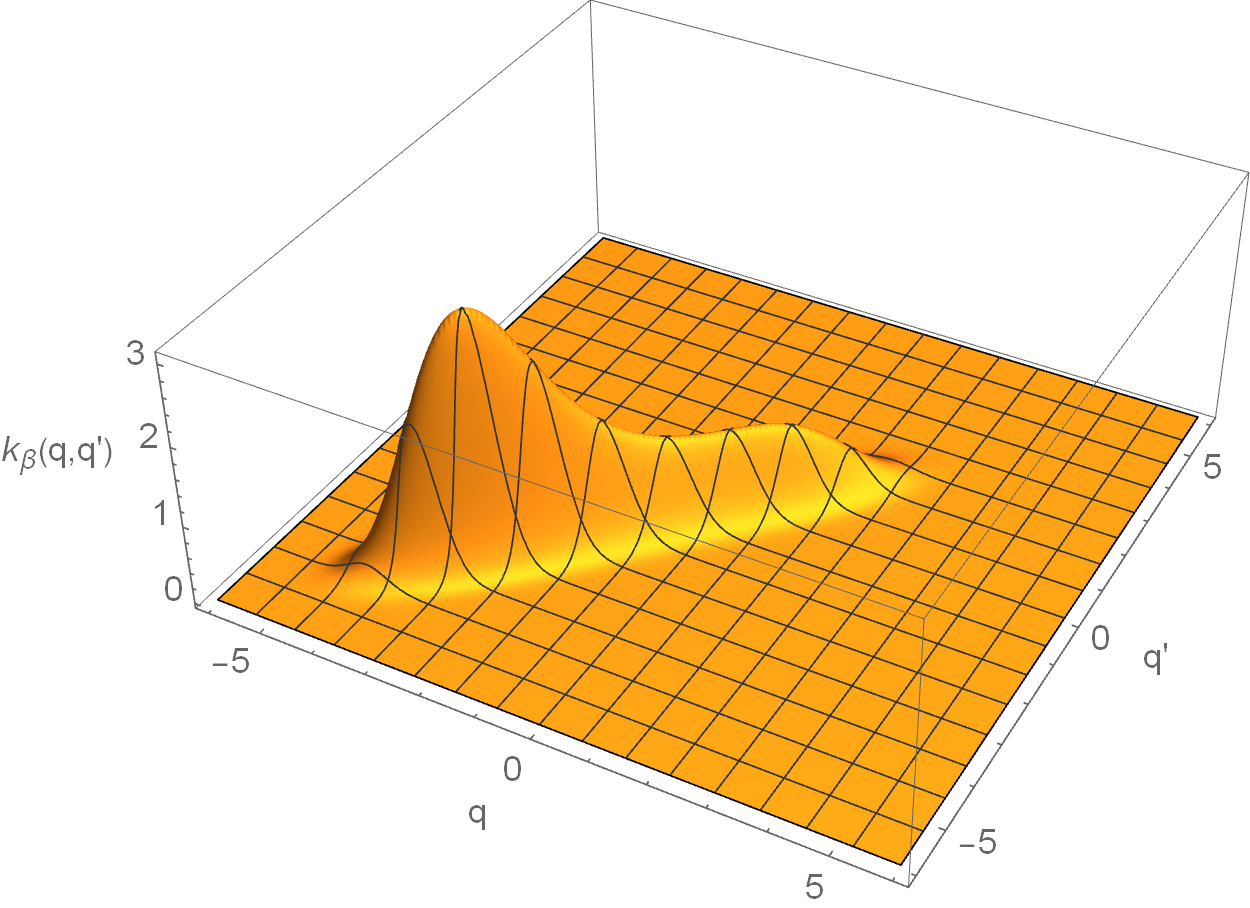}}%
\caption{Integration kernel \eqref{eq:particle_chain_kernel} of the particle chain, for the factorized case $\gamma = 0$ and generic $\gamma = 1$. The remaining parameters are $\lambda = \mu = \frac{1}{5}$ and $\beta = 5$.}
\label{fig:particle_chain_kernel}
\end{figure}

For the following numerical examples, we set $\lambda = \mu = \frac{1}{5}$ and $\eta = 1$. Fig.~\ref{fig:particle_chain_kernel} visualizes the kernel in Eq.~\eqref{eq:particle_chain_kernel}. The factor $\e^{-\frac{1}{2} \beta \gamma (q - q')^2}$ localizes the kernel around the line $q = q'$, which poses a challenge for accurately ``sampling'' it using a limited number of points $z_i$ in Eq.~\eqref{eq:T_beta_def}.

\begin{figure}[!tb]
\centering
\subfloat[free energy]{\includegraphics[width=0.45\textwidth]{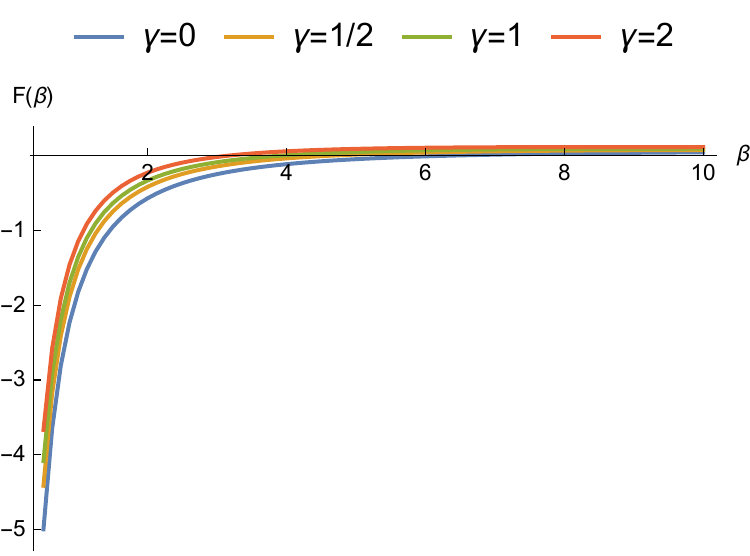}\label{fig:particle_chain_free_energy}}%
\hspace{0.05\textwidth}%
\subfloat[relative error]{\includegraphics[width=0.45\textwidth]{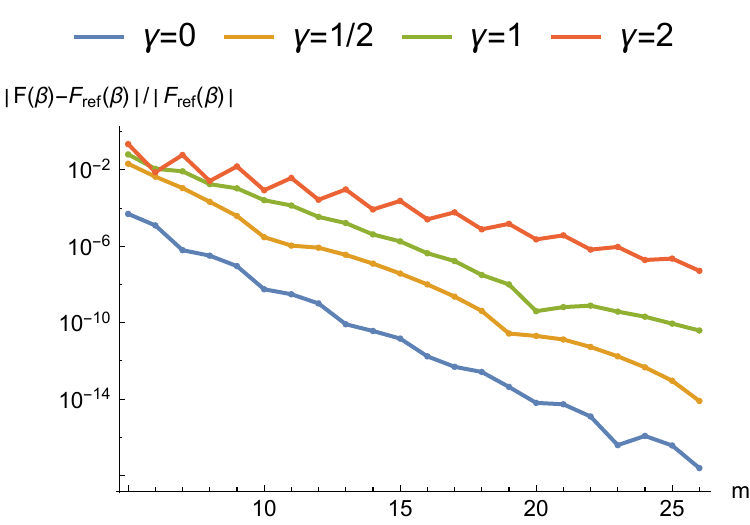}\label{fig:particle_chain_free_energy_convergence}}
\caption{(a) Free energy of the classical particle chain based on a Hamiltonian of the form \eqref{eq:Hpc} with potential \eqref{eq:particle_chain_V}, for parameters $\lambda = \mu = \frac{1}{5}$ and $\eta = 1$. (b) Corresponding convergence plot of the free energy computation in dependence of the number of quadrature points $m$, for $\beta = 5$. The reference values for $\gamma = 0$ have been obtained via Eq.~\eqref{eq:particle_chain_Z1_gamma0}, and for $\gamma > 0$ by using a large $m = 30$.}
\end{figure}

Fig.~\ref{fig:particle_chain_free_energy} shows the free energy as function of $\beta$, for several values of $\gamma$, and Fig.~\ref{fig:particle_chain_free_energy_convergence} the relative error depending on the number of quadrature points $m$. One observes exponential convergence. The less favorable shape of the kernel with increasing $\gamma$, as mentioned above (see also Fig.~\ref{fig:particle_chain_kernel}), translates to a slower convergence rate.

\begin{figure}[!tb]
\centering
\subfloat[$\left\langle \tfrac{1}{2} (q_\ell - q_{\ell+1})^2 \right\rangle$]{\includegraphics[width=0.45\textwidth]{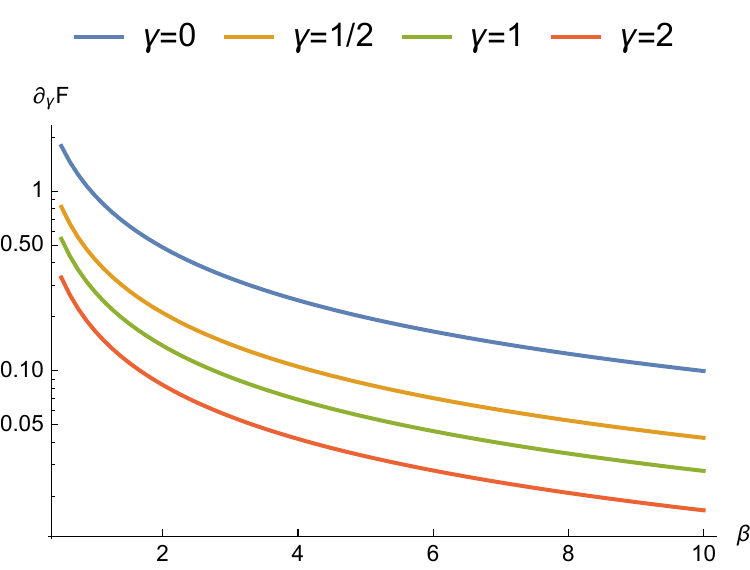}}%
\hspace{0.05\textwidth}%
\subfloat[$\langle e_\ell \rangle$]{\includegraphics[width=0.45\textwidth]{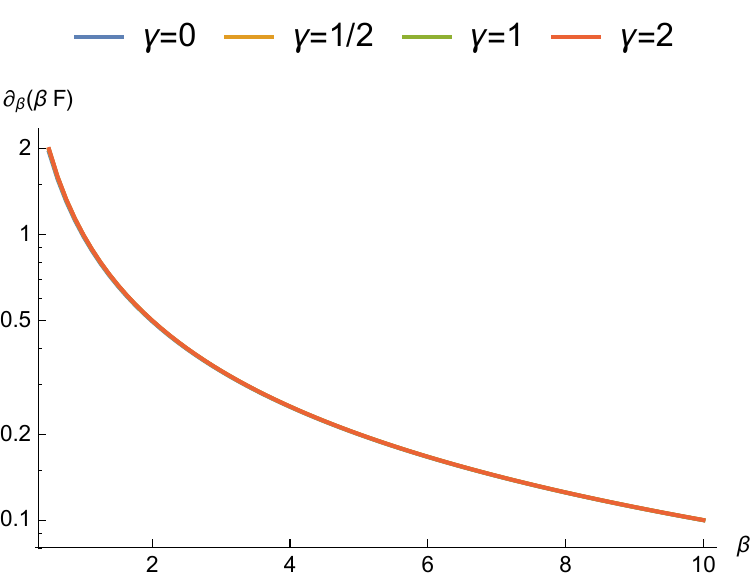}}%
\caption{Average quantities in Eq.~\eqref{eq:particle_chain_avr} based on derivatives of the free energy and plotted on a logarithmic scale in dependence of $\beta$. The derivatives are numerically approximated via a finite difference scheme of order $6$.}
\label{fig:particle_chain_avr}
\end{figure}

Based on the free energy one can obtain averages and higher-order cumulants following the well-known procedure based on derivatives of $F$. For example, the average squared particle distance and energy per site are
\begin{equation}
\label{eq:particle_chain_avr}
\left\langle \tfrac{1}{2} (q_\ell - q_{\ell+1})^2 \right\rangle = \partial_{\gamma} F, \qquad \langle e_\ell \rangle = \partial_{\beta} (\beta F),
\end{equation}
independent of $\ell$ by translation invariance. In practice, a higher-order finite difference scheme on a fine grid is well suited to calculate the derivatives. Fig.~\ref{fig:particle_chain_avr} shows these averages, for the same parameters as before. One notices that the average energy hardly depends on $\gamma$.

As a remark, for the case of vanishing on-site potential, $V_{\text{loc}} \equiv 0$, the model conserves momentum, and the statistical mechanics description changes accordingly \cite{Spohn2014}. Numerically computing the free energy is less challenging in this case since the partition function factorizes after introducing the ``stretch'' $r_{\ell} = q_{\ell+1} - q_\ell$.

\subsection{Discrete nonlinear Schr\"odinger equation}

The method described in Sect.~\ref{sec:method} has been employed in the work \cite{NLS2015} on the discrete nonlinear Schr\"odinger equation. Here we present and elaborate on the numerical aspects in more detail.

To be self-contained, we first restate the physical setup: the central object is a complex-valued wave field $\psi_\ell$ ($\ell = 1, \dots, L$) governed by the semiclassical Hamiltonian
\begin{equation}
\label{eq:NLS_hamiltonian}
H(\psi_1, \dots, \psi_L) = \sum_{\ell=1}^L \big( \tfrac{1}{2} \abs{\psi_{\ell+1} - \psi_\ell}^2 + \tfrac{1}{2}\,g\,\abs{\psi_\ell}^4\big),
\end{equation}
with parameter $g > 0$ (so-called defocusing case). The corresponding partition function reads
\begin{equation}
Z_L(\mu, \beta) = \int \e^{-\beta (H - \mu N)} \, \ud\psi_1 \cdots \ud\psi_L,
\end{equation}
where we have introduced the chemical potential $\mu$ as additional parameter, which is dual to the total particle number $N = \sum_{\ell=1}^L \abs{\psi_\ell}^2$.

\begin{figure}[!tb]
\centering
\includegraphics[width=0.5\textwidth]{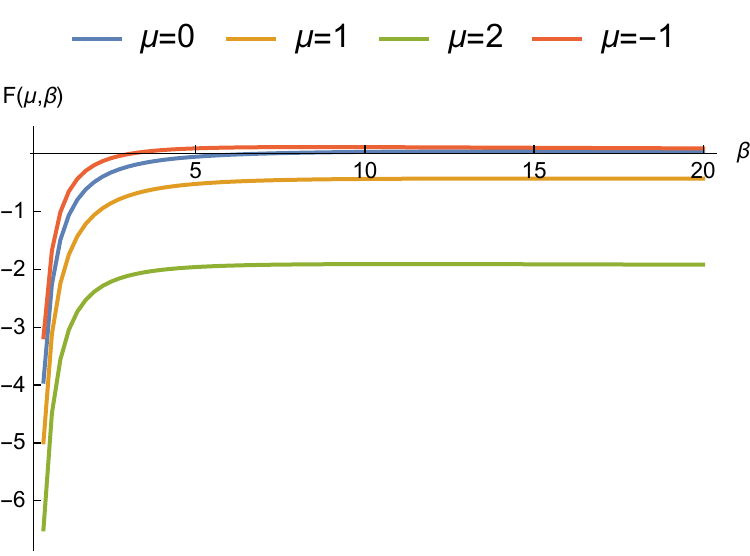}
\caption{Free energy of the discrete nonlinear Schr\"odinger equation governed by the Hamiltonian \eqref{eq:NLS_hamiltonian} on an infinite chain, for parameter $g = 1$.}
\label{fig:NLS_free_energy}
\end{figure}

A symplectic change of variables to polar coordinates leads to the representation
\begin{equation}
\psi_\ell = \sqrt{\rho_\ell}\, \e^{\mathrm{i} \varphi_\ell}
\end{equation}
with $\rho_\ell \in \Omega = [0, \infty)$, and $\varphi_\ell \in S^1$ (unit circle). The Hamiltonian in these variables reads
\begin{align}
\label{eq:polarHamiltonian}
H = \sum_{\ell=1}^L \big( {-} \sqrt{\rho_{\ell+1}\,\rho_\ell}\, \cos(\varphi_{\ell+1} - \varphi_\ell) + \rho_\ell + \tfrac{1}{2}\,g\,\rho_\ell^2 \big).
\end{align}
It depends only on phase differences, which implies the invariance under the global shift $\varphi_\ell \mapsto \varphi_\ell + \phi$. For evaluating the partition function, the $\varphi_\ell$ integrals can be calculated in closed form \cite{RasmussenPRL2000}. This leads to
\begin{equation}
\label{eq:Z_NLS_polar}
Z_L(\mu, \beta) = \e^{\beta \frac{1}{2} \mu^2 L/g} \int_{\Omega} \cdots \int_{\Omega} \prod_{\ell=1}^L k_{\beta}(\rho_\ell, \rho_{\ell+1}) \, \e^{-\beta \frac{1}{2} g \left(\rho_\ell - \frac{\mu}{g}\right)^2} \, \ud\rho_1 \cdots \ud\rho_L
\end{equation}
with the kernel
\begin{equation}
k_{\beta}(\rho, \rho') = 2\pi I_0\big(\beta \sqrt{\rho\,\rho'}\big)\, \e^{-\beta \frac{1}{2} (\rho + \rho')}.
\end{equation}
$I_0$ is the modified Bessel function of the first kind. For this example, we use the second factor of the integrand in Eq.~\eqref{eq:Z_NLS_polar} as weight function:
\begin{equation}
\omega: \Omega \to \R^+, \quad \omega(z) = c \e^{-a (z - b)^2/2}
\end{equation}
with $a = \beta g$, $b = \mu/g$ and the normalization constant
\begin{equation}
\label{eq:weight_normalization_NLS}
c = \frac{2 \sqrt{\frac{a}{2 \pi }}}{1 + \erf(b \sqrt{\frac{a}{2}})}.
\end{equation}
After constructing a Gauss quadrature rule on $\Omega = [0, \infty)$ as described in Sect.~\ref{subsect:nystrom}, we form the symmetric matrix in Eq.~\eqref{eq:T_beta_def}, here denoted $T_{\mu, \beta}$ since it implicitly also depends on $\mu$ via the quadrature points and weights. Then, taking the prefactor in Eq.~\eqref{eq:Z_NLS_polar} and the normalization constant \eqref{eq:weight_normalization_NLS} into account, one arrives at
\begin{equation}
- \beta F(\mu, \beta) = \lim_{L \to \infty} \frac{1}{L} \log Z_L(\mu, \beta) \approx \beta\,\tfrac{1}{2} \tfrac{\mu^2}{g} + \log \lambda_1(T_{\mu, \beta}) - \log c(\mu, \beta).
\end{equation}
Fig.~\ref{fig:NLS_free_energy} shows the free energy as function of $\beta$, for various values of $\mu$.

\begin{figure}[!tb]
\centering
\subfloat[$\beta = 1$]{\includegraphics[width=0.48\textwidth]{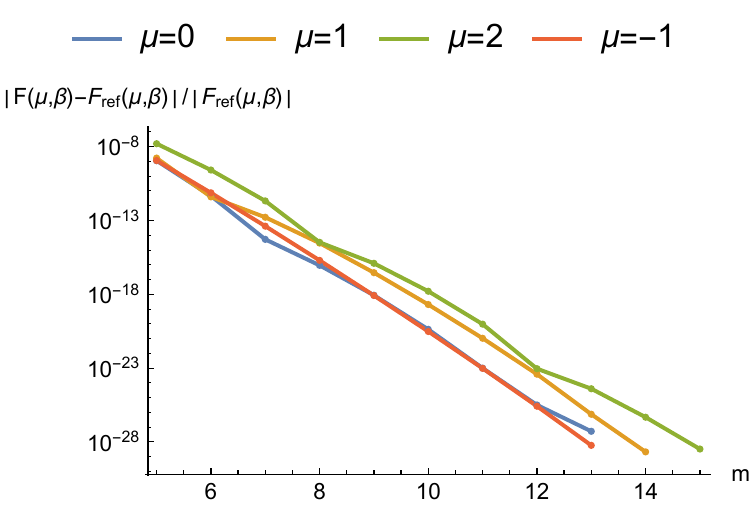}}%
\hspace{0.01\textwidth}%
\subfloat[$\beta = 15$]{\includegraphics[width=0.48\textwidth]{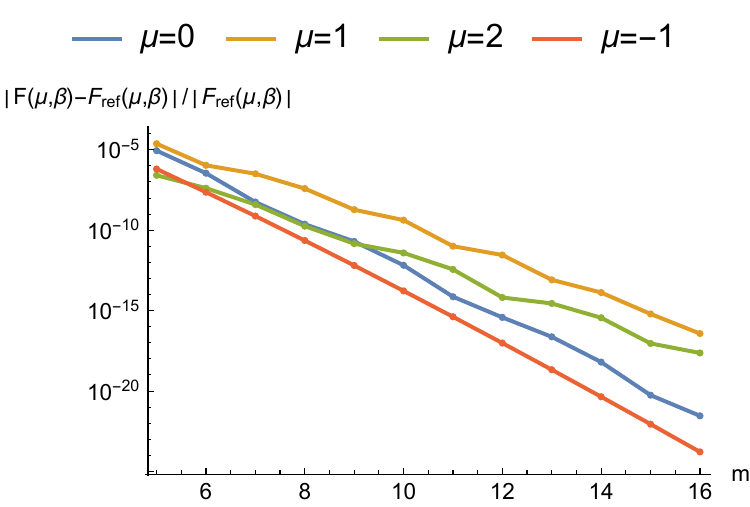}}%
\caption{Relative error of the free energy calculation for the discrete nonlinear Schr\"odinger equation depending on the number of quadrature points $m$, exemplified for $\beta = 1$ and $\beta = 15$.}
\label{fig:NLS_free_energy_convergence}
\end{figure}

Numerically, we again observe exponential convergence with respect to the number of quadrature points, see Fig.~\ref{fig:NLS_free_energy_convergence}. At $\beta = 15$ and $\mu = 1$ for example, $m = 16$ points suffice for double precision accuracy. The reference data stems from a calculation with $m = 20$.

As mentioned before, one can obtain thermodynamic averages based on derivatives of $F(\mu, \beta)$. For example, the average density and energy per lattice site are
\begin{equation}
\langle \rho_\ell \rangle = - \partial_{\mu} F(\mu, \beta), \qquad \langle e_\ell \rangle = \partial_{\beta} (\beta\,F(\mu, \beta)) + \mu\,\langle \rho_\ell \rangle.
\end{equation}
See \cite{NLS2015} for a detailed study of the model.

\subsection{Classical oscillators on a cylindrical lattice}
\label{subsect:lattice}

\begin{figure}[!tb]
\centering
\includegraphics[width=0.6\textwidth]{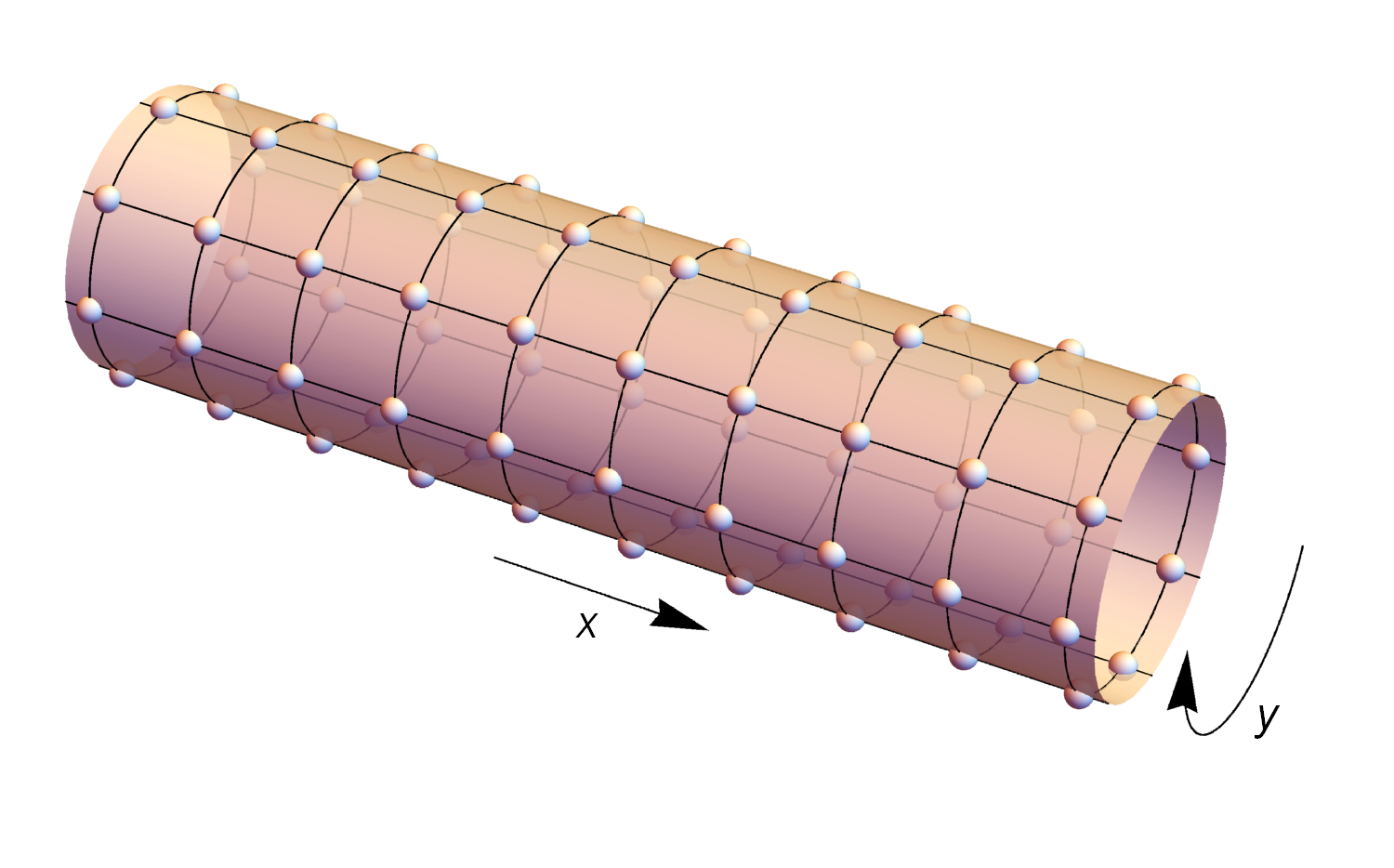}
\caption{Cylindrical lattice topology, using periodic boundary conditions in $y$-direction. Each small dot is a lattice site.}
\label{fig:cylindrical_lattice_topology}
\end{figure}

The numerical method is in principle also applicable to two-dimensional lattices, by using periodic boundary conditions in one direction and reducing the setting to a quasi-one dimensional problem. Specifically, we consider the lattice $\Gamma = \Z/(L_x) \otimes \Z/(L_y)$ for $L_x, L_y \in \N$, i.e., starting with periodic boundary conditions both in $x$- and $y$-direction, but eventually sending $L_x \to \infty$ while keeping $L_y$ finite. Thus we arrive at a cylindrical lattice, as visualized in Fig.~\ref{fig:cylindrical_lattice_topology}.

We identify a lattice site by the index $\ell = (\ell_x, \ell_y) \in \Gamma$, and consider for simplicity scalar spatial variables $q_\ell \in \R$; these could be displacements from the reference positions in one fixed direction, for example. $p_\ell$ denotes the momentum of the $\ell$-th particle.

As demonstration, let the system be governed by the Hamiltonian
\begin{equation}
\label{eq:cylindrical_lattice_hamiltonian}
H = \sum_{\ell \in \Gamma} \left(\tfrac{1}{2} p_\ell^2 + V_{\text{loc}}(q_\ell)\right) + \sum_{\langle \ell, \ell' \rangle} V_{\ell' - \ell}(q_\ell, q_{\ell'}),
\end{equation}
consisting of site-local kinetic and potential energy terms (first sum) as well as nearest neighbor interactions (second sum). Specifically, we consider a local quadratic potential $V_{\text{loc}}(q) = \frac{1}{2} \eta q^2$, $\eta > 0$, and an interaction potential $V_{\Delta\ell}(q, q') = \frac{1}{2} a_{\Delta\ell} (q - q')^2$ with two coefficients $a_{(\pm 1, 0)} = a_x$ and $a_{(0, \pm 1)} = a_y$.

To cast the Hamiltonian \eqref{eq:cylindrical_lattice_hamiltonian} into the form of Eq.~\eqref{eq:formal_hamiltonian}, we subsume the particles contained in one lattice ``ring'' into the vectors
\begin{align}
\vect{p}_{\ell_x} &= \left(p_{(\ell_x, 1)}, \dots, p_{(\ell_x, L_y)}\right) \in \R^{L_y}, \\
\vect{q}_{\ell_x} &= \left(q_{(\ell_x, 1)}, \dots, q_{(\ell_x, L_y)}\right) \in \R^{L_y}
\end{align}
for $\ell_x = 1, \dots, L_x$.

Similar to the particle chain in subsection~\ref{sec:particle_chain}, the momentum integration for evaluating the partition function can be performed explicitly. This results in
\begin{equation}
Z_{(L_x,L_y)}(\beta) = \left(\frac{2\pi}{\beta}\right)^{L_x L_y /2} \tilde{Z}_{(L_x,L_y)}(\beta)
\end{equation}
with
\begin{equation}
\label{eq:cylindrical_lattice_tildeZ_def}
\tilde{Z}_{(L_x,L_y)}(\beta) = \int_{\R^{L_y}} \cdots \int_{\R^{L_y}} \prod_{\ell_x=1}^{L_x} \e^{-\beta \left(\tilde{V}_{\text{loc}}(\vect{q}_{\ell_x}) + \tilde{V}_{\text{int}}(\vect{q}_{\ell_x}, \vect{q}_{\ell_x+1})\right)} \ud^{L_y}q_1 \cdots \ud^{L_y}q_{L_x},
\end{equation}
where
\begin{equation}
\tilde{V}_{\text{loc}}(\vect{q}) = \sum_{\ell_y=1}^{L_y} V_{\text{loc}}(q_{\ell_y})
\end{equation}
and
\begin{equation}
\tilde{V}_{\text{int}}(\vect{q}, \vect{q}') = \sum_{\ell_y=1}^{L_y} \left( \tfrac{1}{2} a_x (q_{\ell_y} - q_{\ell_y}')^2 + \tfrac{1}{4} a_y (q_{\ell_y} - q_{\ell_y+1})^2 + \tfrac{1}{4} a_y (q_{\ell_y}' - q_{\ell_y+1}')^2\right).
\end{equation}
Note that $\tilde{V}_{\text{int}}$ is symmetric in its arguments, i.e., $\tilde{V}_{\text{int}}(\vect{q}, \vect{q}') = \tilde{V}_{\text{int}}(\vect{q}', \vect{q})$ for any $\vect{q}, \vect{q}' \in \R^{L_y}$.

It suggests itself to use the factor $\e^{-\beta \tilde{V}_{\text{loc}}(\vect{q})}$ in Eq.~\eqref{eq:cylindrical_lattice_tildeZ_def} as integration measure, resulting in a tensor product of normal distributions:
\begin{equation}
\label{eq:cylindrical_lattice_weight}
\omega: \R^{L_y} \to \R^+, \quad \omega(\vect{q}) = \frac{\e^{-\frac{1}{2} \beta \eta \norm{\vect{q}}^2}}{\big(\frac{2\pi}{\beta \eta}\big)^{L_y/2}} = \prod_{\ell_y=1}^{L_y} \frac{\e^{-\frac{1}{2} \beta \eta q_{\ell_y}^2}}{\sqrt{2\pi/(\beta \eta)}}.
\end{equation}
We use a rescaled Gauss-Hermite quadrature rule along each coordinate direction, as in Sect.~\ref{sec:particle_chain}. Theorem~\ref{thm:conv} extends straightforwardly to this choice. An alternative, which is less affected by the inherent curse of dimensionality, is a cubature rule dedicated to multidimensional integration \cite{Stroud1971, CoolsRabinowitz1993, Cools1997, Keshavarzzadeh2018}, or sparse grid methods. The convergence properties of such cubature rules are more involved, but they would essentially be inherited by the Nystr\"om method for the dominated eigenvalue. We leave an exploration of these ideas for future work.

Using $\ud\nu(q) = \omega(q)\ud q$ as measure, Eq.~\eqref{eq:cylindrical_lattice_tildeZ_def} becomes
\begin{equation}
\tilde{Z}_{(L_x,L_y)}(\beta) = \left(\frac{2\pi}{\beta \eta}\right)^{L_x L_y/2} \int_{\R^{L_y}} \cdots \int_{\R^{L_y}} \prod_{\ell_x=1}^{L_x} k_{\beta}(\vect{q}_{\ell_x}, \vect{q}_{\ell_x+1}) \, \ud\nu(q_1) \cdots \ud\nu(q_{L_x})
\end{equation}
with the kernel
\begin{equation}
\label{eq:cylindrical_lattice_kernel}
k_{\beta}(\vect{q}, \vect{q}') = \e^{-\beta \tilde{V}_{\text{int}}(\vect{q}, \vect{q}')}.
\end{equation}

Following the factorized quadrature rule, the symmetric matrix in Eq.~\eqref{eq:T_beta_def} takes the form
\begin{equation}
\label{eq:cylindrical_lattice_T_beta}
T_{\beta} = \big(k_{\beta}(\vect{q}_{\vect{i}}, \vect{q}_{\vect{j}})\,\sqrt{w_{\vect{i}}\,w_{\vect{j}}}\big)_{\vect{i},\vect{j}}
\end{equation}
with multi-indices $\vect{i}, \vect{j} \in \{1, \dots, m_0\}^{L_y}$ and the definitions $\vect{q}_{\vect{i}} = (q_{i_1}, \dots, q_{i_{L_y}})$, $\omega_{\vect{i}} = \omega_{i_1} \cdots \omega_{i_{L_y}}$ and $w_i$, $q_i$, $i = 1, \dots, m_0$ the weights and points of the one-dimensional rescaled Gauss-Hermite quadrature rule. Thus the overall number of weights and points is $m = m_0^{L_y}$.

\begin{figure}[!tb]
\centering
\subfloat[free energy]{\includegraphics[width=0.5\textwidth]{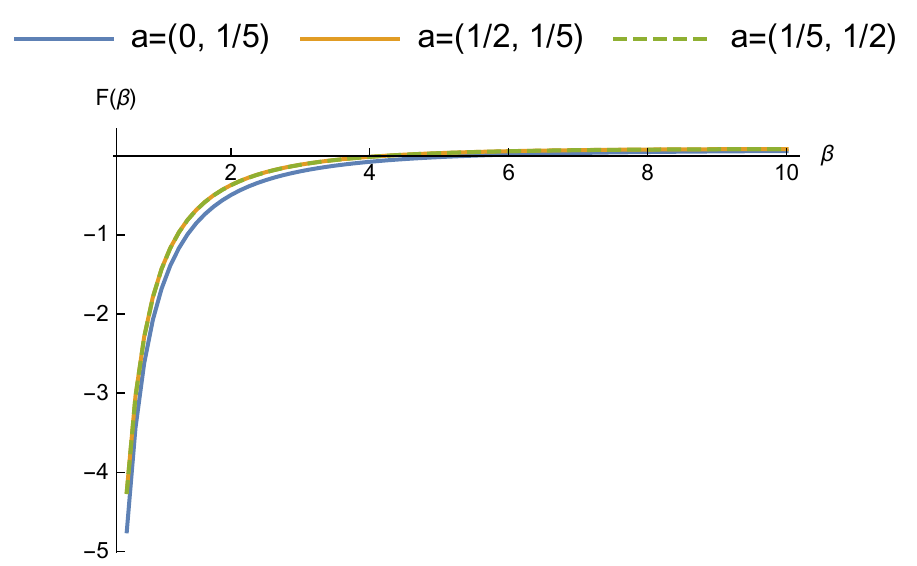}\label{fig:cylindrical_lattice_free_energy}}%
\hspace{0.05\textwidth}%
\subfloat[relative error]{\includegraphics[width=0.45\textwidth]{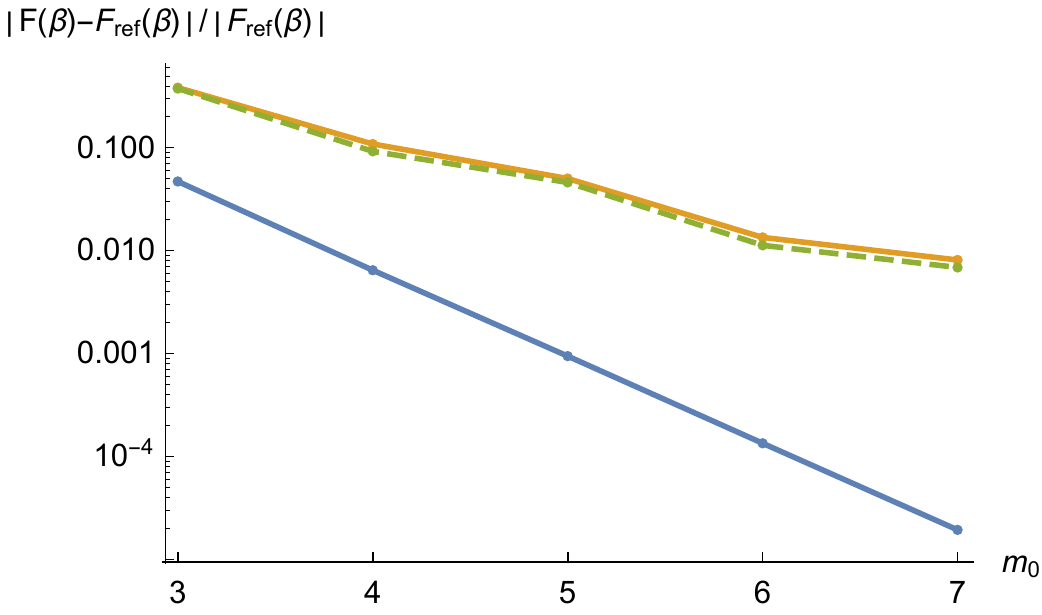}\label{fig:cylindrical_lattice_free_energy_convergence}}
\caption{(a) Free energy of the cylindrical lattice model governed by the Hamiltonian \eqref{eq:cylindrical_lattice_hamiltonian}, for $L_y = 3$ and $\eta = 1$. (b) Corresponding convergence plot of the free energy computation in dependence of the number of quadrature points $m_0$ along one dimension, for $\beta = 5$. The reference values for $a_x = 0$ have been obtained via Eq.~\eqref{eq:cylindrical_lattice_Z1_ax0}, and for $a_x \neq 0$ by using $m_0 = 8$.}
\end{figure}

The numerical approximation of the free energy per lattice site is then
\begin{equation}
- \beta F(\beta) = \lim_{L_x \to \infty} \frac{1}{L_x L_y} \log Z_{(L_x,L_y)}(\beta) \approx \log\frac{2\pi}{\beta} - \frac{1}{2} \log\eta + \frac{1}{L_y} \log \lambda_1(T_{\beta}),
\end{equation}
with $L_y$ kept fixed. For the following examples, we set $L_y = 3$, such that the overall number of quadrature points $m = m_0^{L_y}$ remains manageable up to $m_0 = 8$. Fig.~\ref{fig:cylindrical_lattice_free_energy} visualizes the free energy as function of $\beta$, for several combinations of $a = (a_x, a_y)$. One notices that the curve for $a = (\frac{1}{2}, \frac{1}{5})$ is visually indistinguishable from the case with interchanged parameters $a_x \leftrightarrow a_y$, pointing to the conclusion that the influence of a finite $L_y$ compared to the ``infinite'' $L_x$ on the free energy is quite small.

Fig.~\ref{fig:cylindrical_lattice_free_energy_convergence} shows the corresponding relative error depending on the number of quadrature points $m_0$ along one dimension. In the special case $a_x = 0$ (without coupling in $x$-direction), the partition function factorizes, such that, analogous to Sect.~\ref{sec:particle_chain}, $\tilde{Z}_{(L_x,L_y)}(\beta)\vert_{a_x = 0} = (\tilde{Z}_{(1,L_y)}(\beta)\vert_{a_x = 0})^{L_x}$ with
\begin{equation}
\label{eq:cylindrical_lattice_Z1_ax0}
\tilde{Z}_{(1,L_y)}(\beta)\vert_{a_x = 0} = \int_{\R^{L_y}} \e^{-\beta \left(\frac{1}{2} \eta \norm{\vect{q}}^2 + \sum_{\ell_y=1}^{L_y} \frac{1}{2} a_y (q_{\ell_y} - q_{\ell_y+1})^2 \right)} \ud^{L_y} q.
\end{equation}
We evaluate this integral numerically and use it as reference for computing the relative error in Fig.~\ref{fig:cylindrical_lattice_free_energy_convergence} for $a_x = 0$. The relative error is still rather large for $a = (\frac{1}{2}, \frac{1}{5})$ and $a = (\frac{1}{5}, \frac{1}{2})$; as before, this observation can be explained by the difficulty of accurately sampling the kernel \eqref{eq:cylindrical_lattice_kernel} via \eqref{eq:cylindrical_lattice_T_beta} using a small number of quadrature points along each coordinate. To mitigate this issue for the present example, one could associate the $a_y$ terms in the Hamiltonian to the integration measure $\omega$ instead of the kernel, at the expense of a more complicated quadrature rule.

In summary, this application example demonstrates that our method can in principle handle two-dimensional lattice topologies as well, although the large number of required quadrature points (when interpreting the problem as quasi one-dimensional) limits the size of the periodic dimension $L_y$ in practice.

\section{Conclusion and outlook}

The convergence plots of the example applications illustrate the validity of Theorem~\ref{thm:conv}, which states that approximating the dominant eigenvalue of the discretized kernel inherits the favorable exponential convergence properties of the underlying quadrature rule in the case of kernels extending analytically into the complex domain. 

To optimize the numerical performance of the method for the cylindrical topologies, one could further exploit the factorized structure of $k_{\beta}$ in \eqref{eq:cylindrical_lattice_kernel} along coordinate directions, or use sparse grid methods for the quadrature as mentioned above.

Concerning transfer operator techniques in general, it could be fruitful to adopt ideas from quantum mechanics (see e.g.~\cite{HaegemanVerstraete2017}), or set oriented numerical methods \cite{DellnitzJunge1999}.

\paragraph{Acknowledgments}

C.M.\ likes to thank Herbert Spohn and Abhishek Dhar for helpful discussions, and the Munich Center for Quantum Science and Technology for support.

\section*{Appendix}

We give here an alternative, conceptually much simpler proof\footnote{To the best of our knowledge, this proof has not yet been given in the literature.} of Theorem~\ref{thm:conv} based on the theory of the Fredholm determinant
of the kernel $k_\beta$, namely
\[
D(\mu) = 1 + \sum_{n=1}^\infty \frac{(-\mu)^n}{n!} \int_{\Omega}\cdots\int_{\Omega}\det_{i,j=1}^n k_\beta(z_i,z_j) \,d\nu(z_1)\cdots d\nu(z_n).
\]
Given the assumptions in Sect.~\ref{subsect:assumpt}, $D$ is an entire function whose roots are exactly the reciprocal non-zero eigenvalues of $\mathcal{T}_\beta$, see \cite[\S6.2, Satz 3]{Fen2}. The Weierstrass product \cite[\S6.4, Satz 1]{Fen2}
\[
D(\mu) = e^{\alpha \mu + \beta \mu^2} \prod_n \left(1-\mu \lambda_n(\mathcal{T}_\beta)\right) e^{\mu \lambda_n(\mathcal{T}_\beta)} 
\]
shows that the multiplicities of the roots of $D$ and the multiplicities of the non-zero eigenvalues of $\mathcal{T}_\beta$ agree. In particular, $1/\lambda_1(\mathcal{T}_\beta)$ is a simple root and therefore $D'(\lambda_1(\mathcal{T}_\beta)^{-1})\neq 0$.
Also, \cite[Thm.~6.2]{Bornemann2010} (whose proof can literally be extended to the current assumptions) gives
\begin{equation}\label{eq:error}
\det(I - \mu T_\beta) = D(\mu) + \epsilon_m(\mu)\tag{\#}
\end{equation}
where, uniformly for bounded $\mu$, the error is given by $\epsilon_m(\mu) = O(m^{-k})$ or $\epsilon_m(\mu) = O(e^{-cm})$ according to whether $k_\beta \in C^{k-1,1}$ or $k_\beta$ extends analytically into the complex plane. By Perron--Frobenius $\lambda_1(T_\beta)$ is the simple, dominant, strictly positive eigenvalue of the entry-wise positive matrix $T_\beta$, which by  the argument principle of complex analysis must satisfy
\[
1/\lambda_1(T_\beta) \to 1/\lambda_1(\mathcal{T}_\beta) \qquad (m\to\infty).
\]
Hence, by inserting $\mu=1/\lambda_1(T_\beta)$ into \eqref{eq:error}, followed by a Taylor expansion, we get
\[
\lambda_1(T_\beta) = \lambda_1(\mathcal{T}_\beta) + \frac{\lambda_1(\mathcal{T}_\beta)^2}{D'(\lambda_1(\mathcal{T}_\beta)^{-1})} \epsilon_m(\lambda_1(\mathcal{T}_\beta)^{-1})) + O(\epsilon_m(\lambda_1(\mathcal{T}_\beta)^{-1}))^2),
\]
which completes the proof.

\end{document}